\documentclass[12pt,a4paper]{article}
\usepackage{amssymb,amsfonts,amsmath,amsthm,graphicx}
\usepackage[latin1]{inputenc}

\newenvironment{dem}{\noindent {\bf Proof:}}{\hfill ${\sqcap}\hspace*{-2.8mm}{\sqcup} $ \newline}

\theoremstyle{plain}
\newtheorem{teo}{Theorem}[section]
\newtheorem{defin}[teo]{Definition}
\newtheorem{lema}[teo]{Lemma}
\newtheorem{cor}[teo]{Corollary}
\newtheorem{prop}[teo]{Proposition}

\newtheorem{obs}[teo]{Remark}

\begin{document}

\pagestyle{plain}

\pagenumbering{arabic}

\begin{center}\large{\bf Goldie twisted partial skew power series rings }\footnote{ 
$^{2}$The second named author was partially supported by CAPES}\end{center}

\vskip5mm

\begin{center}{\bf {\rm Wagner $Cortes^1$}}, {\rm Simone  $Ruiz^2$},\end{center}

\begin{center} {\footnotesize $^{1}$ Instituto de Matem\'atica\\
Universidade Federal do Rio Grande do Sul,\\
Porto Alegre-RS, Brazil\\
Av. Bento Gon\c{c}alves, 9500,\\
 91509-900\\
e-mail: {\it wocortes@gmail.com}}
\end{center}

\begin{center}{\footnotesize $^{2}$  Instituto de Matem\'{a}tica\\
Universidade Federal do Rio Grande do Sul\\
91509-900, Porto Alegre-RS, Brazil\\
e-mail:  {\it si-ruiz@hotmail.com}}
\end{center}

\date{}
\begin{abstract}

In this article, we work with unital twisted partial action of $\mathbb{Z}$ on an unital ring $R$ and we introduce the twisted partial skew power series rings and twisted partial skew Laurent series rings.  We study primality,  semi-primality and  prime ideals in these rings.  We completely describe the prime radical  in partial skew Laurent series rings. Moreover, we study the  Goldie property in partial skew power series rings and partial skew Laurent series rings and we describe conditions of the semiprimality of twisted partial skew power series rings.

\end{abstract}

\section*{Introduction}\

Partial actions of groups have been introduced in the theory of
operator algebras as a general approach to study $C^{*}$-algebras by
partial isometries (see, in particular, \cite{E1} and \cite{E3}),
and crossed products classically, as well-pointed out in
\cite{DES1}, are the center of the rich interplay between dynamical
systems and operator algebras (see, for instance, \cite{M1} and
\cite{Q1}). The general notion of (continuous) twisted partial
action of a locally compact group on a $C^{*}$-algebra and the
corresponding crossed product were introduced in \cite{E1}.
Algebraic counterparts for some notions  mentioned above were
introduced and studied in \cite{DE}, stimulating further
investigations, see for instance, \cite{LF},  \cite{CCF}, \cite{laz e mig}
and references therein. In particular, twisted partial actions of
groups on abstract rings and corresponding crossed products were
recently introduced in \cite{DES1}.

 In \cite{wag e fer}, it  was introduced the partial skew polynomial rings and partial skew Laurent of polynomials,  and the authors studied prime and maximal ideals. In  \cite{CFMH}, it  was  investigated the Goldie property in partial skew polynomial rings and partial skew Laurent of polynomial. In \cite{Gobbi}, it was introduced the concept of partial skew power series rings and in the authors studied when it is Bezout and distributive.

 The authors in \cite{Letzter1} and \cite{Letzter2},  studied the Goldie rank  and prime ideals in skew power series ring and skew Laurent series rings with the assumption of noetheriany on the base ring.  In this article, we consider twisted partial actions of $\mathbb{Z}$ and  we introduce the twisted partial skew power series rings and twisted partial skew Laurent series rings $R [[x;\alpha,w]]$ and $R\langle x;\alpha,w\rangle$, respectively, where $\alpha$ is a twisted partial action of $\mathbb{Z}$ on an unital ring $R$.  We  study the Goldie property,  prime ideals, primality and semiprimality  in these rings which generalizes the results presented in \cite{Letzter1} and \cite{Letzter2}.

This article is organized as follows:

In the Section 1, we give some preliminaries and results that will be used during this paper.

In the Section 2,  we study the primality and semiprimality of twisted partial skew power series rings and twisted partial skew Laurent series rings. We describe the prime radical of twisted partial skew Laurent series rings and we study the  prime ideals of these rings.

In the Section 3,  we study the Goldie rank of the twisted partial power series rings and twisted partial skew Laurent series rings  and as a consequence we study the Goldie property of these rings. Morever, we study when the twisted  partial skew power series rings is semiprime and we give a description of the prime radical of twisted partial skew power series rings, when the unital twisted partial action of $\mathbb{Z}$ has enveloping action.

\section{ Preliminaries}

In this section, we recall some notions about twisted partial actions on rings, more details can be found in \cite{DE}, \cite{DES1} and \cite{DES2}. We introduce, in this section, the twisted partial skew power series rings and twisted partial skew Laurent series rings.

From now on,  $R$ will be always an unital ring, unless otherwise stated.

We begin with the following definition that is a particular case of   (\cite{DES2}, Definition 2.1).

\begin{defin} \label{def1}
An unital  \textit{twisted partial action} \index{twisted partial action} of the additive abelian group
$\mathbb{Z}$ on a ring $R$ is a triple
\begin{center}
\vskip-1mm $\alpha = \big(\{D_i\}_{i \in \mathbb{Z}}, \{\alpha_i\}_{i \in \mathbb{Z}}, \{w_{i,j}\}_{(i,j) \in \mathbb{Z}\times \mathbb{Z}}\big)$,
\end{center}
\vskip-1mm where for each $i \in \mathbb{Z}$, $D_i$ is a two-sided ideal in
$R$ generated by a central idempotent $1_i$,  $\alpha_g:D_{-i} \rightarrow D_i$ is an isomorphism of
rings  and for each $(i,j) \in \mathbb{Z} \times \mathbb{Z}$, $w_{i,j}$ is an
invertible element of $D_i D_{i+j}$, satisfying the
following postulates, for all $i,j,k \in \mathbb{Z}$:
\begin{itemize}
\item [$(i)$] \vskip-2.2mm $D_{1} = R$ and $\alpha_{1}$ is the identity map of $R$;
\item [$(ii)$] \vskip-2.2mm $\alpha_i(D_{-i} D_j)= D_i D_{i+j}$;
\item [$(iii)$] \vskip-2.2mm $\alpha _i \circ \alpha _j (a)= w_{i,j} \alpha_{i+j}(a) w_{i,j}^{-1}$, for all $a \in D_{-j} D_{-ji}$;
\item [$(iv)$] \vskip-2.2mm $w_{i,1}=w_{1,i}=1$;
\item [$(v)$] \vskip-2.2mm $\alpha_i(a w _{j,k}) w _{i,j+k}= \alpha_i(a)w_{i,j} w _{i+j,k}$, for all $a \in D_{-i} D_{j}D_{j+k}$.
\end{itemize}
\end{defin}


\begin{obs}

 If $w_{i,j} = 1_i1_{i+j}$, for all $i,j\in \mathbb{Z}$,  then we
have a  partial action  which is a particular case  of
(\cite{DE}, Definition 1.1) and when $D_i = R$, for all $i\in \mathbb{Z}$,
we have that $\alpha$ is a twisted global action.
\end{obs}

Let $\beta = \big(T, \{\beta_i\}_{i \in \mathbb{Z}}, \{u_{i, j}\}_{(i, j)
\in \mathbb{Z} \times \mathbb{Z}}\big)$ be  a twisted global action of a group $\mathbb{Z}$
on a (non-necessarily unital) ring $T$ and $R$ an ideal of $T$
generated by a central idempotent $1_R$. We can restrict $\beta$
to $R$ as follows: putting \linebreak $D_i = R\cap \beta_{i}(R) = R
\beta_{i}(R)$, $i\in \mathbb{Z}$, each $D_{i}$ has an identity element
$1_R\beta_{i}(1_R)$. Then defining
$\alpha_{i}=\beta_{i}|_{D_{-i}}$, $\forall i\in \mathbb{Z}$, the items ($i$),
($ii$) and ($iii$) of \linebreak Definition \ref{def1} are satisfied.
Furthermore, defining $w_{i, j} = u_{i,
j}1_R\beta_i(1_R)\beta_{i+j}(1_R)$, $\forall \,\, i,j\in \mathbb{Z}$, the items ($iv$),
($v$) e ($vi$) of Definition \ref{def1} are also satisfied. So, we obtain a twisted partial
action of $\mathbb{Z}$ on $R$.


The following definition appears in (\cite{DES2}, Definition 2.2).

\begin{defin}\label{def2}
A twisted global action $\big(T, \{\beta_i\}_{i\in \mathbb{Z}},
\{u_{i,j}\}_{(i,j)\in \mathbb{Z}\times \mathbb{Z}}\big)$ of a  group $\mathbb{Z}$ on an
associative (non-necessarily unital) ring $T$  is said to be an
enveloping action \index{enveloping action} (or a globalization)
of an unital  twisted partial action $\alpha$ of $\mathbb{Z}$ on a ring $R$ if, 
there exists a monomorphism $\varphi:R\rightarrow T$ such that,
for all $i$ and $j$ in $\mathbb{Z}$:
\begin{itemize}
\item [$(i)$] \vskip-2.2mm $\varphi(R)$  is an  ideal of $T$;
\item [$(ii)$] \vskip-2.2mm $T = \displaystyle\sum_{i\in \mathbb{Z}}\beta_i(\varphi(R))$;
\item [$(iii)$] \vskip-2.2mm $\varphi(D_i) = \varphi(R)\cap \beta_i(\varphi(R))$;
\item [$(iv)$] \vskip-2.2mm $\varphi\circ \alpha_i(a) = \beta_i\circ \varphi(a)$, for all $a\in D_{-i}$;
\item [$(v)$] \vskip-2.2mm $\varphi(aw_{i,j}) = \varphi(a)u_{i,j}$ and $\varphi(w_{i,j}a) = u_{i,j}\varphi(a)$, for all $a\in D_iD_{i+j}$.
\end{itemize}
\end{defin}

In (\cite{DES2}, Theorem 4.1), the authors studied  necessary and
sufficient conditions for an unital  twisted partial action $\alpha$ of a
group $\mathbb{Z}$ on a ring $R$ has an enveloping action. Moreover, they
studied which rings satisfy such conditions.

Suppose that $(R, \alpha, w)$ has an enveloping action
$(T,\beta,u)$. In this case, we may assume that $R$ is an ideal of
$T$ and we can rewrite the conditions of the Definition \ref{def2}
as follows:
\begin{itemize}
\item [$(i')$] \vskip-3mm $R$  is an  ideal of $T$;
\item [$(ii')$] \vskip-3mm $T = \displaystyle\sum_{i\in \mathbb{Z}}\beta_i(R)$;
\item [$(iii')$] \vskip-3mm $D_i = R\cap \beta_i(R)$, for all $i\in \mathbb{Z}$;
\item [$(iv')$] \vskip-3mm $\alpha_i(a) = \beta_i(a)$, for all $x\in D_{-i}$ and  $i\in \mathbb{Z}$;
\item [$(v')$] \vskip-3mm $aw_{i,j} = au_{i,j}$ and $w_{i,j}a = u_{i,j}a$, for all $a\in D_iD_{i+j}$ and $i, j \in \mathbb{Z}$.
\end{itemize}

Given an unital  twisted partial action $\alpha$ of $\mathbb{Z}$  on a ring
$R$,  we define  the twisted partial skew Laurent series  rings
$R\langle x;\alpha,w\rangle=\displaystyle\bigoplus_{i\in \mathbb{Z}} D_ix^i$ whose elements are the series
\begin{center}
\vskip-1mm $\displaystyle\sum_{j\geq s} a_jx^j$, with $a_j\in D_j$
\end{center}
\vskip-1mm
 with the usual addition and multiplication defined by
\vskip-5mm
$$(a_ix^i)(a_jx^j) = \alpha_i(\alpha^{-1}_i(a_i)b_j)w_{i,j}x^{i+j}.$$
\vskip2mm

Using the similar techiniques of  (\cite{DES1}, Theorem 2.4), $R\langle x;\alpha,w\rangle$ is an associative
ring whose identity is $1_Rx^{0}$. Note that, we have the
injective morphism $\phi: R\rightarrow R\langle x;\alpha,w\rangle$, defined
by $r\mapsto rx^0$ and we can consider $R\langle x;\alpha,w\rangle$ as
an extension of $R$.  Moreover, we consider the twisted partial power series rings as    a subring of  $R\langle x;\alpha,w\rangle$ which we denote it by $R [[x;\alpha,w]]$   whose elements are the series $\displaystyle\sum_{i\geq 0} b_ix^i$ with sum and multiplication rule defined as before.



Let $\alpha$ be an unital  twisted partial action of a  group $Z$ on a ring
$R$. An ideal  $S$ of $R$ is said to be $\alpha$-ideal ($\alpha$-invariant ideal) if,
$\alpha_i(S\cap D_{-i}) \subseteq S\cap D_i$, for all $i\geq
0$  $(\alpha_i(S\cap D_{-i}) = S\cap D_i$, for all $i\in
\mathbb{Z})$.

If $S$ is an $\alpha$-ideal ($\alpha$-invariant ideal), then we have the ideals
\begin{center}
$S [[x;\alpha,w]]=\left\{\displaystyle\sum_{i\geq 0}a_ix^i\, |\, a_i\in S\cap D_i\right\}$ ($S\langle x;\alpha,w\rangle=\left\{\displaystyle\sum_{i\geq m}a_ix^i\, |\, a_i\in S\cap D_i \,\, m\in\mathbb{Z}\right\}$)
\end{center}
is an ideal of  $R [[x;\alpha,w]]$ $(R\langle x;\alpha,w\rangle)$. Note that,  if $I$ is a right  ideal of $R$, then $I [[x;\alpha,w]]=\{\displaystyle\sum_{i\geq 0} a_ix^i: a_i\in D_i\}$ and $I\langle x;\alpha,w\rangle=\{\displaystyle\sum_{i\geq m}b_ix^i:b_i\in D_i\}$ are right ideals of $R [[x;\alpha,w]]$ and $R\langle x;\alpha,w\rangle$, respectively.

Note that for each  $\alpha$-invariant ideal $I$ of $R$,  the unital twisted partial action $\alpha$ can be extended to an unital  twisted partial action $\overline{\alpha}$ of $\mathbb{Z}$ on $R/I$ as follows: for each $i\in \mathbb{Z}$, we define $\overline{\alpha}_i:D_{-i}+I\longrightarrow D_i+I$, putting $\overline{\alpha}_i(a+I)=\alpha_{i}(a)+I$, for all $a\in D_{-i}$, and for each $(i,j)\in \mathbb{Z}\times \mathbb{Z}$, we extend each $w_{i,j}$ to $R/I$ by  $\overline{w}_{i,j}=w_{i,j}+I$.


Moreover, when  $(R,\alpha,w)$ has enveloping action $(T,\beta,u)$, then by similar methods presented in  Section 2 of \cite{laz e mig}, we have that $(T/I^{e}, \overline{\beta},\overline{u})$ is the enveloping action of $(R/I,\overline{\alpha},\overline{w})$, where $I^e$ is the $\beta$-invariant ideal such that  $I^{e}\cap R=I$.

We finish this section with some comments about twisted partial actions of finite type that will be necessary in this paper.

The following definition is a particular case of  (\cite{lmsw}, Definition 4.13).

\begin{defin} Let $\alpha$ be an unital  twisted partial action. We say that $\alpha$  is of finite type if, there exists a finite subset $\{s_1,s_2, \cdots ,s_n\}$ of $\mathbb{Z}$ such that \begin{center} $\displaystyle\sum_{i=1}^n D_{j+s_i} =R$,\end{center} for all $j\in \mathbb{Z}$. \end{defin}

 It is convenient to point out that in the same way as in (\cite{laz e mig}, Proposition 1.2)  as proved  in \cite{lmsw}, we have  that an unital  twisted partial action $\alpha$ of  $\mathbb{Z}$ on an unital ring $R$ with an enveloping action $(T,\beta,u)$ is of finite type if, and only if,  there exists $s_1,\cdots, s_n\in \mathbb{Z}$ such that $T=\displaystyle\sum_{i=1}^n \beta_{s_{i}} (R)$  and this is equivalent to say that $T$  has an identity element.


\section{Primality and semiprimality}

In this section, $\alpha$ will denote an unital twisted partial action of $\mathbb{Z}$ on an unital ring $R$, unless otherwise stated.
We begin this section with the following  proposition, whose proof is standard,  and we put it here  for the sake of completeness.  





\begin{prop}\label{quociente} If  $I$ is an  $\alpha$-invariant ideal of  $R$, then   $\frac{R [[x;\alpha,w]]}{I [[x;\alpha,w]]}\simeq (\frac{R}{I}) [[x;\overline{\alpha},\overline{w}]]$. Moreover, the same result holds to  $R\langle x;\alpha,w\rangle$.
\end{prop}
\begin{dem} 



   We define   $\varphi:\frac{R [[x;\alpha,w]]}{I [[x;\alpha,w]]}\rightarrow  (\frac{R}{I}) [[x;\overline{\alpha},\overline{w}]]$  by $\varphi(\displaystyle\sum_{i\geq 0}a_ix^i+I[[x;\alpha,w]])=\displaystyle\sum_{i\geq 0} (a_i+ I)x^i$. We easily have that $\varphi$ is an isomorphism.  So, $\frac{R [[x;\alpha,w]]}{I [[x;\alpha,w]]}\simeq (\frac{R}{I}) [[x;\overline{\alpha},\overline{w}]]$
\end{dem}

The following definition firstly appeared in \cite{wag e fer}  for ordinary partial actions

\begin{defin} Let $\alpha$ be an unital  twisted partial action of $\mathbb{Z}$ on $R$ and $I$ an ideal of $R$.

(i) $I$ is $\alpha$-prime if, $I$ is an $\alpha$-invariant ideal and for each $J$ and $K$ $\alpha$-invariant ideals of $R$ such that $JK\subseteq I$ implies that either $J\subseteq I$ or $K\subseteq I$.

(ii) $I$ is strongly $\alpha$-prime if,  $I$ is $\alpha$-invariant and  for each ideal $M$ of $R$  and $\alpha$-ideal $N$ of $R$ such that $MN\subseteq I$ implies that either $M\subseteq I$ or $N\subseteq I$.

\end{defin}


Let $a\in R$. Then  we define the $\alpha$-invariant ideal generated by $a$  as  $J=\displaystyle\sum_{i\in \mathbb{Z}}R\alpha_i(a1_{-i})R$.

In the next result, we study necessary and sufficient conditions for $\alpha$-primality and strongly $\alpha$-primality.

\begin{lema}\label {primo2}  (1) Let $P$ be an $\alpha$-invariant ideal of $R$.  The following conditions are equivalent:

(a)    $P$ is $\alpha$-prime 

(b)  For each  $a, b\in R$ such that  $\alpha_j(a1_{-j})R\alpha_i(b1_{-i})\subseteq P$, for all  $i, j\in \mathbb{Z}$, then either  $a\in P$ or  $b\in P$

(c) $R/P$ is $\overline{\alpha}-prime$, where $\overline{\alpha}$ is the extension of twisted partial action $\alpha$ to $R/P$

(2) Let $P$ be an $\alpha$-invariant ideal of $R$.  The following conditions are  equivalent:

(a) $P$ is strongly $\alpha$-prime 

(b)  For each  $a, b\in R$ such that  $aR\alpha_j(b1_{-j})\subseteq P$, for all  $j\geq 0$, then either  $a\in P$ or  $b\in P$.

(c) $R/P$ is strongly $\overline{\alpha}$-prime, where $\overline{\alpha}$ is the extension of twisted partial $\alpha$ to $R/P$. \end{lema}

\begin{dem} (1)   $(a)\Rightarrow (b)$

Let  $a, b\in R$ such that  $\alpha_j(a1_{-j})R\alpha_i(b1_{-i})\subseteq P$, for all  $i, j\in \mathbb{Z}$. Then, if we fix  $j$ we have that  $$\alpha_j(a1_{-j})\displaystyle\sum_{i\in \mathbb{Z}}R\alpha_i(b1_{-i})R\subseteq P.$$  and consequently,   we get  $$\displaystyle\sum_{j\in \mathbb{Z}}R\alpha_j(a1_{-j})R\displaystyle\sum_{i\in \mathbb{Z}}R\alpha_i(b1_{-i})R\subseteq P.$$ Since the ideals  $\displaystyle\sum_{i\in \mathbb{Z}}R\alpha_j(a1_{-j})R$ and  $\sum_{i\in \mathbb{Z}}R\alpha_i(b1_{-i})R$ are  $\alpha$-invariant, then, by assumption, we have that  either $\displaystyle\sum_{i\in \mathbb{Z}} R\alpha_j(a1_{-j})R\subseteq P$ or  $\displaystyle\sum_{i\in \mathbb{Z}}R\alpha_i(b1_{-i})R\subseteq P$. So, either $a\in P$ or  $b\in P$.

$(b)\Rightarrow (a)$

Let $I,J$ be $\alpha$-invariant ideals of $R$ such that $IJ\subseteq P$, take $a\in I$ and suppose that there exists $b\in J\setminus P$. Then, $(\displaystyle\sum_{i\in \mathbb{Z}}R\alpha_i(a1_{-i})R)(\displaystyle\sum_{j\in \mathbb{Z}}R\alpha_j(b1_{-j})R)\subseteq P$. Thus, by assumption, we have that either \begin{center} $\displaystyle\sum_{i\in \mathbb{Z}}R\alpha_i(a1_{-i})R\subseteq P$  or $\displaystyle\sum_{j\in \mathbb{Z}}R\alpha_j(b1_{-j})R\subseteq P$. \end{center}  Hence, $a\in P$, because $b\notin P$.  So, $I\subseteq P$. 

$(a)\Rightarrow (c)$ 

Let  $a, b\in R$ such that  \begin{center} $\overline{\alpha}_j((a+P)(1_{-j}+P))(R/P)\overline{\alpha}_i((b+P)(1_{-i}+P))=\overline{0}$, \end{center} for all  $i, j\in \mathbb{Z}$. Then,   $\alpha_j(a1_{-j})R\alpha_i(b1_{-i})\subseteq P$, for all  $i, j\in \mathbb{Z}$. Thus, by assumption, we have that either $a\in P$ or $b\in P$. So, either $a+P=\overline{0}$ or $b+P=\overline{0}$. 

$(c)\Rightarrow (a)$ 

Let $I$ and $J$ be $\alpha$-invariant ideals of $R$ such that $IJ\subseteq P$. Thus, $\overline{I}\overline{J}=\overline{0}$ in $R/P$. Hence, by assumption, we have that either $\overline{I}=\overline{0}$ or $\overline{J}=\overline{0}$. So, either $I\subseteq P$ or $J\subseteq P$.

The proof of item (ii) is analogous.
      \end{dem}

It is convenient to point out that $R$ is $\alpha$-prime  (strongly $\alpha$-prime) if the zero ideal is $\alpha$-prime (strongly $\alpha$-prime). Next, we have an easy consequence of Lemma \ref{primo2}.

\begin{lema}\label{anelprimo} (i) Let $\alpha$ be an unital twisted partial action of $\mathbb{Z}$ on $R$. Then $R$ is  $\alpha$-prime if, and only if, for each  $a, b\in R$ such that  $\alpha_j(a1_{-j})R\alpha_i(b1_{-i})=0$ for all  $i,j\in \mathbb{Z}$, we have that either  $a=0$ or  $b=0$.

(ii) Let $\alpha$ be an unital twisted partial action of $\mathbb{Z}$ on $R$. Then, $R$ is  strongly $\alpha$-prime if, and only if, for each  $a, b\in R$ such that  $aR\alpha_i(b1_{-i})=0$ for all  $i\geq 0$, we have that either  $a=0$ or  $b=0$
\end{lema}


It is convenient to point out that if, $L$ is a nonzero right ideal of $R\langle x;\alpha,w\rangle$, then $L\cap R [[x;\alpha,w]]$ is a nonzero right ideal of $R [[x;\alpha,w]]$ because of for each nonzero element $f\in L$, there exists $s\geq 0$ such that $0\neq f1_sx^s \in L\cap R [[x;\alpha,w]]$. Moreover,  if  a right ideal $M$ of $R\langle x;\alpha,w\rangle$  is such that $M\cap R [[x;\alpha,w]]=0$, then we have that $M=0$.   We use these facts without further mention.

In the next result, we study conditions for the primality  of $R [[x;\alpha,w]]$ and $R\langle x;\alpha,w\rangle$ which partially generalizes (\cite{Letzter2}, Propositions 2.5 and 2.7).

  \begin{prop}\label{anelprimo}The following statements hold.
\begin{description}
              \item[(a)] $R$ is  $\alpha$-prime  if and only if $R\langle x;\alpha,w\rangle$ is prime.
              \item[(b)]  $R [[x;\alpha,w]]$ is  prime if and only if $R$ is  strongly $\alpha$-prime. In particular, if $R [[x;\alpha,w]]$ is prime, then $R$ is $\alpha$-prime.
              \item[(c)] If  $R [[x;\alpha,w]]$ is prime, then  $R\langle x;\alpha,w\rangle$ is prime

              \end{description}
              \end{prop}

              \begin{dem}
              \begin{description}
\item[(a)] Suppose that $R\langle x;\alpha,w\rangle$ is prime and let   $I$ and  $J$ be $\alpha$-invariant ideals of  $R$ such that  $IJ=0$. Then  $$I\langle x;  \alpha,w\rangle J\langle x;\alpha,w\rangle \subseteq (IJ)\langle x;\alpha,w\rangle =0.$$ By the fact that   $R\langle x;\alpha,w\rangle$ is  prime, we have that either  $I\langle x;\alpha,w\rangle=0$ or  $J\langle x;\alpha,w\rangle=0$.  Hence, either $I=0$ or $J=0$. So,  $R$ is  $\alpha$-prime.

    Conversely, let  $f, g \in R\langle x;\alpha,w\rangle$ be nonzero elements,   suppose that  $fR\langle x;\alpha,w \rangle g=0$ and consider  $m$ and  $n$ the smallest  integers such that  $f_m\neq 0$ and  $g_n\neq 0$ where  $f=\displaystyle\sum_{i\geq m}f_ix^i$ and  $g=\displaystyle\sum_{i\geq n}g_ix^i$.  Note that, for each   $i\in \mathbb{Z}$, $fD_ix^ig\subseteq f R\langle x;\alpha,w\rangle g=0$ and  we have that $f_mD_i\alpha_i(1_{-i}g_n)=0$,  for all  $i\in \mathbb{Z}$. Hence,  for each  $j\in \mathbb{Z}$, we have that  $\alpha_{j}(f_m1_{-j})R\alpha_i(g_n1_{-i})=0$,
 for all  $i\in \mathbb{Z}$.

Consequently,  by  Lemma \eqref{anelprimo}, we have that  $f_m=0$ or $g_n=0$, which is a contradiction. So,  $R\langle x;\alpha,w\rangle$ is prime.

\item[(b)] The proof is similar of the item (a).

\item[(c)] Let  $I$ and  $J$  be ideals of  $R\langle x;\alpha,w\rangle$ such that  $IJ=0$. Thus,  $$ (I\cap R [[x;\alpha,w]])(J\cap R [[x;\alpha,w]])=0.$$
Since  $I\cap R [[x;\alpha,w]]$ and  $J\cap R [[x;\alpha,w]]$ are ideals of  $R [[x;\alpha,w]]$, then   we have that  either $I\cap R [[x;\alpha,w]]=0$ or  $J\cap R [[x;\alpha,w]]=0$. Hence either  $I=0$ or $J=0$. So,  $R\langle x;\alpha,w\rangle$ is prime.

\end{description}
\end{dem}

\begin{obs} The authors in (\cite{Letzter2}, Propositions 2.5 and 2.7)  used noetherianity to get the  equivalences mentioned there. For us  to obtain the same equivalences  in the case of twisted partial actions, we would need to know  if the following question has  a  positive answer, but until now we do not know.

\begin{center} Are all $\alpha$-ideals of $R$ $\alpha$-invariant ideals when $R$ is Noetherian?\end{center}

So, if this question has a positive answer we would have that  $R$ is $\alpha$-prime $\Leftrightarrow$ $R [[x;\alpha,w]]$ is prime $\Leftrightarrow$ $R\langle x;\alpha,w\rangle$ is prime.

\end{obs}

The following result is a direct consequence  of the last proposition.

              \begin{cor}  If  $R$ is a prime ring, then $R [[x;\alpha,w]]$ is a prime ring.\end{cor}

             The proof of the following result is similar to Proposition \ref{anelprimo} and it partially generalizes (\cite{Letzter2},  Corollary 2.12)

\begin{cor} \label{primeideals} The following statements hold.

(a) Suppose that $I$ is $\alpha$-invariant ideal of $R$. Then  $I$ is  $\alpha$-prime if and only if  $I\langle x;\alpha,w\rangle$ is prime.

(b) Suppose that $I$ is an $\alpha$-invariant ideal of $R$. Then $I$ is strongly $\alpha$-prime if and only if $I [[x;\alpha,w]]$ is prime.

\end{cor}

The following result generalizes  (\cite{Letzter1}, Theorem 3.18) and  is a direct consequence of the last corollary.

\begin{cor} Let $\alpha$ be an unital  twisted partial action of  $\mathbb{Z}$ on $R$ and $I$ a strongly  $\alpha$-prime ideal of $R$. Then, there exists a prime ideal $P$ of $R [[x;\alpha,w]]$ such that $P\cap R=I$. Moreover, if $I$ is $\alpha$-prime, then there exists a prime ideal $Q$ of $R\langle x;\alpha,w\rangle$ such that $Q\cap R=I$. \end{cor}




In (\cite{Letzter2}, Proposition 2.11) is used the noetherianity property to prove the result in the case of skew Laurent series rings, but in that proof the assumption was not necessary.   The next result  generalizes (\cite{Letzter2}, Proposition 2.11).

\begin{prop} \label{primality1}  If  $K$ is a prime ideal of  $R\langle x;\alpha,w\rangle$, then  $K\cap R$ is an $\alpha$-prime  ideal of $R$.

\end{prop}

\begin{dem} Let $K$ be an prime ideal  of $R\langle x;\alpha,w\rangle$. Then, we easily have that $K\cap R$ is an ideal of $R$. We claim that   $K\cap R$ is an $\alpha$-prime ideal of $R$. In fact, let $a\in (K\cap R)\cap D_{-i}$, for $i\in \mathbb{Z}$. Then,  $1_ix^iax^{-i}\in K$.   Thus,
$$1_ix^iax^{-i}=1_i\alpha_i(a)w_{i, -i}=\alpha_i(a)w_{i, -i} \in K\cap D_i$$ and since $w_{i, -i}$ is an invertible element of $D_i$, we get that  $\alpha_i(a)w_{i,-i}w_{i,-i}^{-1} \in (K\cap R)\cap D_i)$. Hence,   $\alpha_i(a)\in (K\cap R)\cap D_i$ and it follows that  $\alpha_i((K\cap R)\cap  D_{-i})\subseteq (K\cap R)\cap D_i$. By similar methods, we show that  $\alpha_{i}^{-1}((K\cap R)\cap  D_i)\subseteq (K\cap R)\cap D_{-i}$. Consequently,   $\alpha_i((K\cap R)\cap D_{-i})= (K\cap R)\cap D_i$, for all  $i\in \mathbb{Z}$ and we have that  $K\cap R$ is an $\alpha$-invariant ideal  of  $R$.

By Proposition 2.1 we have that \begin{center} $\Psi:(R/(K\cap R))\langle x;\overline{\alpha},\overline{w}\rangle \rightarrow (R\langle x;\alpha,w\rangle)/((K\cap R)\langle x; \alpha,w\rangle)$\end{center}  defined by $\Psi(\sum_{i\geq s} \overline{a_i}x^ i)=\sum_{i\geq s}a_ix^i+(K\cap R)\langle x;\alpha,w\rangle$ is an isomorphism. Note that $K/((K\cap R)\langle x;\alpha,w\rangle)$ is a prime ideal and we have that $\Psi^{-1}(K/((K\cap R)\langle x;\alpha,w\rangle))=\overline{K}=\{\sum_{i\geq s} (a_i+(K\cap R))x^i:\sum_{i\geq s} a_ix^i\in K\}$ is a prime ideal in  $(R/(K\cap R))\langle x;\overline{\alpha},\overline{w}\rangle$ and $\overline{K}\cap (R/(K\cap R))=\overline{0}$. Thus, we may   assume that $K\cap R =0$ and in this case we only need to show that $R$ is $\alpha$-prime. In fact, let $I$ and $J$ be $\alpha$-invariant ideals of $R$ such that $IJ=0$.  Hence, $IR\langle x;\alpha,w\rangle J\langle x;\alpha,w \rangle\subseteq I\langle x;\alpha,w\rangle J\langle x;\alpha,w\rangle \subseteq (IJ)\langle x;\alpha,w\rangle=0\subseteq K$. By the fact that $K$ is a prime ideal we have that either $IR\langle x;\alpha,w\rangle\subseteq K$ or  $J\langle x;\alpha,w \rangle\subseteq K$ and it follows that  either $I\subseteq K$ or $J\subseteq K$. So, either $I=0$ or $J=0$ and we have that $R$ is $\alpha$-prime.

\end{dem}

The following notion appears in \cite{CM}.

\begin{defin} Let $\alpha$ be an unital  twisted partial action of $\mathbb{Z}$ on $R$. Then the $\alpha$-nil radical $N_{\alpha}(R)$ of $R$ is the intersection of all $\alpha$-prime ideals of $R$.\end{defin}

From now on, for a ring $S$ we denote its  prime radical by $Nil_{*}(S)$.

Now, we are in conditions to describe the prime radical of $R\langle x;\alpha,w\rangle$.

\begin{prop}\label{primeradical5} Let $\alpha$ be an unital  twisted partial action of $\mathbb{Z}$ on $R$. Then $Nil_{*}(R\langle x;\alpha,w\rangle)=Nil_{\alpha}(R)\langle x;\alpha,w\rangle$. \end{prop}

\begin{dem} Let $P$ be a prime ideal of $R\langle x;\alpha,w\rangle$. Then,  by Proposition \ref{primality1},  we have that $P\cap R$ is $\alpha$-prime. Thus, $Nil_{*}(R\langle x;\alpha,w\rangle)\supseteq Nil_{\alpha}(R)\langle x;\alpha,w\rangle$.

On the other hand, let $I$ be an $\alpha$-prime ideal of $R$. Then, by Corollary \ref{primeideals}, we have that $I\langle x;\alpha,w\rangle$ is prime. Hence, $Nil_{\alpha}(R)\langle x;\alpha,w\rangle\supseteq Nil_{*}(R\langle x;\alpha,w\rangle)$. So,  $Nil_{*}(R\langle x;\alpha,w\rangle=Nil_{\alpha}(R)\langle x;\alpha,w\rangle$. \end{dem}

\begin{prop} \label{semiprime1}  Let $\alpha$ be an unital  twisted  partial action  of $\mathbb{Z}$ on $R$.

(i) If  $R$ is semiprime, then $R\langle x;\alpha,w\rangle$ is semiprime. Moreover, if $R$ is Noetherian and $R\langle x;\alpha,w\rangle$ is semiprime, then $R$ is semiprime.

(ii) Let $I$ be an $\alpha$-invariant ideal of $R$. If $I$ is semiprime, then $I\langle x;\alpha,w\rangle$ is semiprime.

\end{prop}

\begin{dem}
(i)  Assume, by the way of contradiction, that there exists $f=\displaystyle\sum_{i\geq s}f_ix^i$    such that $fR\langle x;\alpha,w\rangle f=0$, where  $f_s\neq 0$. Take any $c\in D_{s}$ and write $b=\alpha_{-s}(c)$ , for some $b\in D_{-s}$. Thus, $fbx^{-s}f=0$ and we have that  $f_s\alpha_{s}^{-1}(b)w_{s,-s}f_s = 0$ . Hence, $f_scw_{-s,s}f_s = 0$ and we get that $f_sD_sf_s = 0$. Since $R$ is a semiprime ring, then $D_s$  is also a semiprime ring. Consequently, $f_s=0$ because $f_s\in  D_s$, a contradiction.  So, $R\langle x;\alpha,w\rangle$ is semiprime.

For the second part, since $R$ is Noetherian, then by  (\cite{lam1}, Theorem 4.10.30)  the prime radical $Nil_{*}(R)$ is nilpotent. As a consequence, there exists $n\geq 1$ such that for every $\alpha$-prime ideal $P$ of $R$ we have that $Nil_{*}(R)^n\subseteq P$ and it follows that  $Nil_{*}(R)\subseteq P$, for every $\alpha$-prime ideal of $R$, because of   (\cite{laz e mig}, Remark 3.2) says that $Nil_{*}(R)$ is an $\alpha$-invariant ideal of $R$.  Hence,   we get that  $Nil_{*}(R)\subseteq Nil_{\alpha}(R)$. By assumption and Proposition \ref{primeradical5}  we have that $Nil_{\alpha}(R)=0$ and consequently, $Nil_{*}(R)=0$ So, $R$ is semiprime.

(ii) The proof is similar of the item (i).

\end{dem}

Fron now on,  we  proceed to give a more close description of the prime ideals of
$R[[x;\alpha,w]]$ and $R\langle x;\alpha,w\rangle$.  The proof of the next result is similar to (\cite{wag e fer},  Proposition 2.6).

\begin{prop}  \label{lem25} Let $P$ be a prime ideal of
$R[[x;\alpha,w]]$ (resp. $R\langle x;\alpha,w\rangle$). Then we have one of the
following possibilities:

(i) $P=Q\oplus \displaystyle\sum_{i\geq 1} D_ix^{i}$, where $Q$ is a prime
ideal of $R$

(resp. $P=Q\oplus \displaystyle\sum_{i\not =0} D_ix^{i}$, where $Q$ is a prime
ideal of $R$ with $D_j\subseteq Q$, for any $j\not =0$).

(ii) $1_ix^{i}\notin P$, for some $i\geq 1$.
\end{prop}

It is clear that for any prime ideal $Q$ of $R$, the ideal $Q\oplus
\displaystyle\sum_{i\geq 1} D_ix^{i}$ is a prime ideal of $R[[x;\alpha,w]]$. Thus,
we are in the case (i) of Proposition \ref{lem25} . If, in addition,
$D_j\subseteq Q$, for all $j\not =0$, it is easy to see that
$P=Q\oplus \displaystyle\sum_{i\not =0} D_ix^{i}$ is an ideal of $R\langle x;\alpha,w\rangle$
which is obviously prime.

From now on, we proceed to study the case  of the item (ii) of the last proposition and we have the following two results.

\begin{prop} Let $P$ be an ideal of $R\langle x;\alpha,w\rangle$.  If  $P\cap R$ is $\alpha$-prime and either
$P=(P\cap R)\langle x;\alpha,w\rangle$ or  $P$ is maximal amongst the ideals $N$ of
$R\langle x;\alpha,w\rangle$ with $N\cap R=P\cap R$, then $P$ is prime.
\end{prop}

\begin{dem}  If $P=(P\cap R)\langle x;\alpha,w\rangle$, then the result follows from Corollary \ref{primeideals}. Now, suppose that $P\neq (P\cap R)\langle x;\alpha,w\rangle$ and let $I,J$ be ideals of $R\langle x;\alpha,w\rangle$ such that $IJ\subseteq P$.  Suppose that  $I\nsubseteq P$ and $J\nsubseteq P$ and we get that  $P\subsetneq I+P$ and $P\subsetneq J+P$. Note that $((I+P)\cap R)((J+P)\cap R)\subseteq P\cap R$ because of $(I+P)(J+P)\subseteq P$. By assumption, we have that either $((I+P)\cap R)\subseteq P\cap R$ or $((J+P)\cap R)\subseteq P\cap R$. Thus, either $(I+P)\cap R=P\cap R$ or $(J+P)\cap R=P\cap R$, which  contradicts the assumption on $P$. Hence,  either $I\subseteq P$ or $J\subseteq P$. So, $P$ is prime.
\end{dem}

The proof of the following  result is similar to the proof of the last proposition.

\begin{prop}  Let $P$ be an ideal of $R[[x;\alpha,w]]$ such that
$1_ix^{i}\notin P$, for some $i\geq 1$ and $P\cap R$ is an $\alpha$-invariant ideal.  If $P\cap R$ is prime and either $P=(P\cap
R)[[x;\alpha,w]]$ or $P$ is maximal amongst the ideals $N$ of
$R[[x;\alpha,w]]$ with $N\cap R=P\cap R$, then $P$ is prime.
\end{prop}

We finish this section with the following remark.

\begin{obs}  Until now, we do not know if it is true or not  the following natural converse of the last two propositions:

(i) If $P$ is a prime ideal of $R\langle x;\alpha,w\rangle$ and $P\neq (P\cap R)\langle x;\alpha,w\rangle$, then $P$ is maximal amongst the ideals $N$ of
$R\langle x;\alpha,w\rangle$ with $N\cap R=P\cap R$.

(ii)   Let $P$ be an ideal of $R[[x;\alpha,w]]$ such that
$1_ix^{i}\notin P$, for some $i\geq 1$ and $P\cap R$ is a strongly $\alpha$-prime ideal of $R$.  If $P$ is a prime ideal of $R [[x;\alpha,w]]$ and $P\neq (P\cap R) [[x;\alpha,w]]$,  then $P$ is maximal amongst the ideals $N$ of
$R[[x;\alpha,w]]$ with $N\cap R=P\cap R$.

.

\end{obs}

\section{Goldie twisted partial skew power series rings}

In this section, $\alpha$ is an unital twisted partial action of $\mathbb{Z}$ on $R$, unless otherwise stated.

Let $S$ be a ring and $M$ a right  $S$-module. We remind that $M$ is uniform if, the intersection of any two nonzero submodules is nonzero,  see (\cite{mcconnel robson}, pg.  52) for more details. According to (\cite{mcconnel robson}, pg. 57)   a ring $S$ is right Goldie if satisfies ACC on right annihilator ideals  and $S$ does not have an infinite direct sum of right uniform   ideals.  In this section, we study the Goldie property in twisted  partial skew Laurent series rings and twisted partial skew power series rings. We begin with the following lemma that will be important to prove the principal results of this section, which generalizes (\cite{Letzter2}, Lemma 2.8).


\begin{lema}\label{submoduloordenado} Let  $V$ be  a right simple  $R$-module. Then  $VR [[x;\alpha,w]]$ is a right $R$-module whose the only   submodules are ordered  in the form  $$ VR[[x; \alpha,w]]\supset V(\sum_{i\geq 1}D_ix^i)\supset V(\sum_{i\geq 2}D_ix^i)\supset \ldots.$$
\end{lema}
\begin{dem} We easily have that    $VR [[x;\alpha,w]]$ is a  right $R [[x;\alpha,w]]$-module   and  note that $ VR [[x;\alpha,w]]\supset V(\displaystyle\sum_{i\geq 1}D_ix^i)\supset V(\displaystyle\sum_{i\geq 2} D_ix^i)\supset \ldots$.

Let  $S$ be a  $R [[x;\alpha,w]]$- submodule of  $VR [[x;\alpha,w]]$ such that $S \neq V\displaystyle\sum_{i\geq 1}D_ix^i$ and  $f=\displaystyle\sum_{i\geq 0} v_i x^i$ a nonzero element of  $S$ with   $0\neq v_0\in V$. Since  $V$ is a simple right  $R$-module, then  $v_0R=V$. Thus there exists  $a_i\in R$ such that  $v_i=v_0a_i$ for all  $i\geq 1$. Let  $g=1+ u_1x+ u_2x^2+ \ldots$ be an element of  $R [[x;\alpha,w]]$ such that
\begin{eqnarray*}fg&=& (v_0 + v_0a_1x + v_0a_2x^2+ \ldots)(1 + u_1x + u_2x^2+ \ldots)\\
&=& v_0 + (v_0u_1 + v_0a_1)x + (v_0u_2+ \alpha_1(\alpha^{-1}_1(v_0a_1)u_1)w_{1,1}+ v_0a_2)x^2 \\
&+& (v_0u_3 + \alpha_1(\alpha^{-1}_1(v_0a_1)u_2)w_{1,2} + \alpha_2(\alpha^{-1}_2(v_0a_2)u_1)w_{2,1} + v_0a_3)x^3+ \ldots.
\end{eqnarray*}
If we take   $u_1=-a_1$, $u_2=-a_2-a_1\alpha_1(u_11_{-1})w_{1,1}$, $u_3=-a_3-a_2\alpha_2(u_11_{-2})w_{2,1}-a_1\alpha_1(u_21_{-1})w_{1,2}$, ..., $u_n= a_n-a_{n-1}\alpha_{n-1}(u_11_{-n+1})w_{n-1,1}-...-a_1\alpha_1(u_{n-1}1_{-1})w_{1,n-1}$.., then we get that   $fg=v_0$.

Note that  $$VR [[x;\alpha,w]]=(v_0R)R [[x;\alpha,w]]\subseteq v_0R [[x;\alpha,w]]=fgR [[x;\alpha,w]]\subseteq fR [[x;\alpha,w]]$$ and  $fR [[x;\alpha,w]]\subseteq VR [[x;\alpha,w]]$. Hence, $VR [[x;\alpha,w]]=fR [[x;\alpha,w]]$, for all  $f\in S$ and it follows that  $VR [[x;\alpha,w]]\subseteq SR [[x;\alpha,w]]\subseteq S$. So,  $V\displaystyle\sum_{i\geq 1}D_ix^i$ is the unique submodule of  $VR [[x;\alpha,w]]$. Finally, following this technique we get the result.
\end{dem}

 Next, we study the uniformity of $VR [[x;\alpha,w]]$ and $VR\langle x;\alpha,w\rangle$.

\begin{prop}\label{idealsimples} Suppose that  $V$ is a right simple ideal of  $R$. The following statements hold.
\begin{description}
 \item[(a)]$VR [[x;\alpha,w]]$ is uniform as  $R [[x;\alpha,w]]$-module.
 \item[(b)]$VR\langle x; \alpha,w\rangle$ is uniform as $R\langle x;\alpha,w\rangle$-module.
\end{description}
\end{prop}
\begin{dem}
\begin{description}
\item[(a)] By Lemma  \eqref{submoduloordenado}, the unique submodules of  $VR [[x;\alpha,w]]$  are $V\displaystyle\sum_{i\geq m}D_ix^i$,  for $m\geq 0$ and note that   $V\displaystyle\sum_{j\geq i}D_jx^j\supset V\displaystyle\sum_{j\geq i+s}D_jx^j$, for all  $s\geq 0$. Thus, $$V\displaystyle\sum_{j\geq s}D_jx^j \cap V\displaystyle\sum_{j\geq t}D_jx^j=V\displaystyle\sum_{j\geq t}D_jx^j\neq 0,$$ always that   $s\geq t$. So,  $VR [[x;\alpha,w]]$ is uniform.

\item[(b)] Let  $L$ be a nonzero submodule of  $VR\langle x;\alpha,w\rangle$. Then,  $L\cap VR [[x;\alpha,w]]$ is a nonzero submodule of  $VR [[x;\alpha,w]]$.  Thus, for each nonzero submodules  $C$ and  $D$ of  $VR\langle x;\alpha,w\rangle$, we have that  $C\cap VR [[x;\alpha,w]]\neq 0$ and  $D\cap VR [[x;\alpha,w]]\neq 0$, and it follows that  $(C\cap D)\cap VR [[x;\alpha,w]]=(C\cap VR [[x;\alpha,w]])\cap (D\cap VR [[x;\alpha,w]])\neq 0$. Hence, $C\cap D\neq0$. So,  $ VR\langle x;\alpha,w\rangle$ is uniform.
\end{description}
\end{dem}

According to (\cite{mcconnel robson}, 2.2.10)  the right Goldie rank of a ring  $S$ is $n$ if there exists a direct sum $\displaystyle\bigoplus_{i=1}^n  I_i$ of uniform right submodules of $S$ such that  $\displaystyle\bigoplus_{i=1}^n  I_i$ is right essential in $S$ and we denote it by $rankS$.

 In (\cite{Letzter2}, Theorem 2.8)  the authors used the noetherianity to prove it.  In  next result we replace the noetherianity condition for a weaker condition, that is, Goldie property  and it generalizes (\cite{Letzter2}, Theorem 2.8).

\begin{teo}\label{dimensaouniforme} If  $R$ is semiprime Goldie, then  $rankR=rankR [[x;\alpha,w]]=rank R\langle x;\alpha,w\rangle$
\end{teo}
\begin{dem} By the fact that  $R$ is semiprime Goldie we have, by (\cite{mcconnel robson}, Theorem 2.3.6) , that  there exists the classical quotient  ring  $E$ of $R$ which is semisimple.  Note that   $rankR=rankE$, because of (\cite{mcconnel robson}. Lemma 2.2.12).    Since  $R\subseteq R\langle x;\alpha,w\rangle\subseteq E\langle x;\alpha^{*},w^{*}\rangle$, then   $$rank E=rank R\leq rank R\langle x;\alpha,w\rangle \leq rank E\langle x;\alpha^{*},w^{*}\rangle,$$ where $\alpha^{*}$ is the extension of the unital twisted partial action $\alpha$ of $R$ to $E$, see (\cite{lmsw}, Theorem 3.12).
Let  $d=rankR$ and  we may suppose without loss of generality that  $R=E$ and $\alpha=\alpha^{*}$.  Then, we can write  $$R=V_1\oplus \cdots \oplus V_d$$ where  $V_i$ is a simple right ideal of  $R$, for all  $i=1, \ldots, d$.

Hence, $$R\langle x;\alpha,w\rangle=V_1R\langle x;\alpha,w\rangle\oplus \cdots \oplus V_dR\langle x;\alpha,w\rangle.$$ and by  Proposition  \eqref{idealsimples}, item $(b)$,  each  $V_iR\langle x;\alpha,w\rangle$ is uniform as right  $R\langle x;\alpha,w\rangle$- module. So,  $rank R\langle x;\alpha,w\rangle=d$.

By similar methods, we have that  $R [[x;\alpha,w]]= V_1R [[x;\alpha,w]]\oplus \cdots \oplus  V_dR[[x, \alpha.w]]$ and by Proposition  \eqref{idealsimples} item $(b)$,  each  $V_iR [[x;\alpha,w]]$ is an unifom submodule of  $R [[x;\alpha,w]]$, for all  $i=1, \ldots d$. So, $rank R [[x;\alpha,w]]=d$.

\end{dem}






Let $S$ be a ring and $a\in S$. The right annihilator of $a$ in $S$ is $Ann_S(a)=\{x\in S: ax=0\}$. Moreover, according to (\cite{mcconnel robson}, Definition 2.2.4) the singular ideal of $S$ is $Z(S)=\{a\in S: Ann_S(a)\,\, is \,\, right \,\, essential\,\, in \,\, S\}$,

Now, we are ready to prove the second principal result of this section.

\begin{teo}\label{goldie} Let $R$ be a semiprime ring. The following conditions are equivalent:
\begin{description}
\item[(a)] $R$ is  Goldie.
\item[(b)] $R [[x;\alpha,w]]$ is  Goldie.
\item[(c)] $R\langle x;\alpha,w\rangle$ is  Goldie.
\end{description}
\end{teo}
\begin{dem}
$(a)\Rightarrow (c)$ By assumption,   Theorem  \ref{dimensaouniforme} and by  Proposition \ref{semiprime1}, item $(i)$  we have that $rankR\langle x;\alpha,w\rangle=rankR<\infty$ and $R\langle x;\alpha,w\rangle$ is semiprime.  We claim that  $R\langle x;\alpha,w\rangle$ is nonsingular. In fact, let  $f\in Z(R\langle x;\alpha,w\rangle)$,  where $f= a_{-j}x^{-j}+ \ldots + a_0+ a_1x+ \ldots$ and
 $I$ a nonzero right ideal of  $R$. Then  $I\langle x;\alpha,w\rangle$ is a right ideal of  $R\langle x;\alpha,w\rangle$  and we obtain that  $Ann_{R\langle x;\alpha,w\rangle}(f)\cap I\langle x;\alpha,w\rangle\neq 0$.  Thus, there exists  $0\neq h\in I\langle x;\alpha,w\rangle\cap  Ann_{R\langle x;\alpha,w\rangle}(f)$, i.e., $fh=0$. We consider,  $h=b_{-k}x^{-k}+ \ldots + b_0 + b_1x+ \ldots $ and suppose without loss of generality that  $b_{-k}\neq 0$. Hence, looking at  the smallest degree of the product  $fh$  we get
$$ a_{-j}\alpha_{-j}(1_jb_{-k})w_{-j, -k} x^{-j-k}=0,$$ which implies that  $a_{-j}\alpha_{-j}(1_jb_{-k})=0$. Consequently, $$\alpha^{-1}_{-j}(a_{-j})\alpha^{-1}_{-j}(\alpha_{-j}(1_jb_{-k}))=0 \Longrightarrow \alpha^{-1}_{-j}(a_{-j})1_jb_{-k}=0 $$ and we have that  $\alpha^{-1}_{-j}(a_{-j})b_{-k}=0$. So, $0\neq b_{-k}\in Ann_R(\alpha^{-1}_{-j}(a_{-j}))$ and we obtain that   $Ann_R(\alpha^{-1}_{-j}(a_{-j}))\cap I\neq 0$ which concludes that  $\alpha^{-1}_{-j}(a_{-j})\in Z(R)$.   By the fact that $R$ is Goldie we have that $\alpha^{-1}_{-j}(a_{-j})=0$. Since,   $\alpha^{-1}_{-j}$ is an isomorphism, then  $a_{-j}=0$. Now, following the similar methods,  we obtain  that    $f=0$. Hence,  $Z(R\langle x;\alpha,w\rangle)=0$.  Therefore, by (\cite{mcconnel robson}, Theorem 2.3.6)  we get that   $R\langle x;\alpha,w\rangle$ is Goldie.

We need to show that  $(c)\Rightarrow (b)\Rightarrow (a)$. In fact, note that  $$R\subset R [[x;\alpha,w]]\subset R\langle x;\alpha,w\rangle$$  and by the fact that $R$ is semiprime and Goldie,  we have, by Theorem \ref{dimensaouniforme}, that  $rank R=rankR [[x;\alpha,w]]=rankR\langle x;\alpha,w\rangle.$   Since the chain conditions on right annihilators is inherited by subrings we obtain the desired result.
\end{dem}

In the article
\cite{lmsw},  the authors worked with twisted partial actions of finite type and the rings  satisfied some finiteness conditions as Goldie property. But,  at that time the authors did not notice such assumption would imply the existence of the enveloping action. So, in the  next result, we show that the unital twisted partial actions  on algebras with finite Goldie rank that are of finite type,  have enveloping action.

\begin{teo}\label{envolvente}  Let $R$ be a ring with finite uniform dimension and $\alpha$ a twisted partial action of $\mathbb{Z}$ on $R$.  If $\alpha$ is of finite type, then $\alpha$ has enveloping action.
\end{teo}

\begin{dem}  By assumption, there exists a finite set  $\{g_1, \ldots , g_n\}$ of  $\mathbb{Z}$ such that  $$R=D_{g+g_1} + \ldots + D_{g+g_n},$$
for every  $g\in \mathbb{Z}$.
We claim that  $R$ can be written as direct sum of indecomposable rings. In fact,   each $D_{g_i}$ has identity $1_{g_i}$ and by similar methods of (\cite{laz e mig}, Remark 1.11)  we can write  $$R=F_1\oplus\ldots \oplus F_n,$$  where each  $F_i$ is an ideal of  $D_{g+g_i}$, $i=1,\ldots, n$,  generated by a central idempotent. Now, if each  $F_i$ is indecomposable we are done. Next, if there exists $1\leq j\leq n$ such that $F_j$ is not indecomposable,  then we may write  $F_j=F_j^{1}\oplus  F_j^{2}$,  and we get  $$R= F_1\oplus...\oplus F_j^{1}\oplus F_j^{2}\oplus... F_n.$$
Proceeding in this manner with all other decomposable components we may write $$R=A_1\oplus ...\oplus A_n$$ Now if all $A_i$   are indecomposable, then  we are done. If it is not,  proceed with similar methods as before. Since $rankR$ is finite, then  the process must stop and we have that $R$ is a direct sum of indecomposable rings where, each one is generated by a central idempotent of $R$. So, by (\cite{DES2}, Theorem 7.2), $(R,\alpha,w)$ has enveloping action.

\end{dem}

Let $\alpha$ be an unital  twisted partial action of $\mathbb{Z}$ on $R$ that admits enveloping action $(T,\beta,u)$. Following \cite {DE} and \cite{DES2},  we exhibit an explicit Morita context between
$R\langle x;\alpha,w\rangle$ and $T\langle x; \beta,u \rangle$ whose restriction to $T[[x;\beta,u]]$
gives also a Morita context between $R[[x;\alpha,w]]$ and
$T[[x;\beta,u]]$.

Recall that given two rings $R$ and $S$, bimodules ${_R}U_S$ and
${_S}V_R$ and maps $\theta:U\otimes_{S} V\rightarrow R$ and
$\psi:V\otimes_{R} U\rightarrow S$, the collection
$(R,S,U,V,\theta,\psi)$ is said to be a Morita context if the
array
\[ \left[ \begin{array}{cc}
R & V \\
U & S
\end{array} \right], \]
with the usual formal operations of $2\times 2$ matrices, is a
ring.

The following result is proved in (\cite{mcconnel robson}, Theorem 3.6.2), for
rings with identity element. Actually, in the proof of the result, it
is not used the fact that the rings have identity element and the
modules $U$ and $V$ are unital modules. So, we can easily see that
the following is true for rings which do not necessarily have
identity.

\begin{teo} \label{morita} Let $(R,S,U,V,\theta,\psi)$ be a Morita context.
Then there is an order preserving one-to-one correspondence
between the sets of prime ideals $P$ of $R$ with $P\nsupseteqq UV$
and prime ideals $P^{'}$ of $S$ with $P^{'}\nsupseteqq VU$. The
correspondence is given by $P\longmapsto \{s\in S:UsV\subseteq
P\}$ and $P^{'}\longmapsto \{r\in R: VrU\subseteq
P^{'}\}$.\end{teo}

Following the similar ideas of (\cite{DE}, Section 5), we put $U=\{\displaystyle\sum_{i\in
\mathbb{Z}}a_{i}x^{i}: \, \, a_{i}\in R\, \, ,{\rm for}\, \, {\rm
all}\,\, i\in \mathbb{Z}\}$ and $V=\{\displaystyle\sum_{i\in
\mathbb{Z}}a_{i}x^{i}: a_{i}\in \beta_{i}(R)\, \, ,{\rm for}\, \,
{\rm all}\, \, i\in \mathbb{Z}\}$. Then, it can easily be seen that
$UT\langle x;\beta,u \rangle \subseteq U$, $T\langle x;\beta,u\rangle V\subseteq V$,
$R\langle x;\alpha,w\rangle U\subseteq U$ and $VR\langle x;\alpha,w\rangle\subseteq V$ (to show
the relations recall that $\beta_{j}(R)$ is an ideal of $T$ and
$S_{j}=\beta_{j}(S_{-j})$). In case we want to consider
$R[[x;\alpha,w]]$ and $T[[x;\beta,u]]$ we restrict $U$ and $V$ to have
just power series  and we have similar relations.

Thus, we have the Morita contexts
$(R\langle x;\alpha,w\rangle,T\langle x;\beta,u\rangle ,U,V,\theta,\psi)$ and
$(R[[x;\alpha,w]],T[[x;\beta,u]]],U,V,\theta,\psi)$, where $\theta$ and
$\psi$ are obvious.

The proof of the following lemma is similar to (\cite{wag e fer}, Lemma 2.2).

\begin{lema} \label{prim1} Let $P$ be a prime ideal of  $R[[x;\alpha,w]]$.
Then, there exists a unique prime ideal $P^{'}$ of  $T[[x;\beta,u]]$, given by Theorem \ref{morita},  which satisfies $P^{'}\cap R[[x;\alpha,w]]]=P$.
\end{lema}

We have the following easy consequence.

\begin{cor} \label{prim2}  There is a one-to-one correspondence, via contraction,
between the set of all prime ideals of $R[[x;\alpha,w]]$ and the set
of all prime ideals of $T[[x;\beta,u]]$ which do not contain $R$.
\end{cor}

The next result is important to prove the last main result of this article and it is an easy consequence  of Lemma \ref{prim1} and Corollary \ref{prim2}.

\begin{cor}\label{primeradical} Let $\alpha$ be an unital  twisted partial action of $\mathbb{Z}$ on $R$ with enveloping action $(T,\beta,u)$. Then $Nil_{*}(T [[x;\beta,u]])\cap R [[x;\alpha,w]]=Nil_{*}(R [[x;\alpha,w]])$.  \end{cor}

 Based on the last results,  we will proceed to describe  the prime radical of $R [[x;\alpha,w]]$ when $(R,\alpha,w)$ has enveloping action $(T,\beta,u)$ and for this we need the following result.



\begin{lema}\label{primeradical1} Let $\beta$ be a twisted global action of $\mathbb{Z}$ on a ring $S$ with cocycle $u$. Then,  the prime radical $Nil_*(S[[x;\beta,u]])$ of  $S[[x,\beta, u]]$ is $Nil_*(S[[x;\beta,u]])=Nil_{*}(S)\cap N_{\beta}(S)\oplus \displaystyle\sum_{i\geq 1}N_{\beta}(S)x^i$, where $N_{\beta}(S)$ is the intersection of all strongly $\beta$-prime ideals  of $S$.\end{lema}

\begin{dem} We have two classes of prime ideals in $S [[x;\beta,u]]$, i.e., \begin{center} $\mathcal{F}_1=\{ P: \,\, prime \,\, ideal\,\, such \,\, that\,\ , S [[x;\beta,u]]x\subseteq P\} $ \end{center} and \begin{center}$\mathcal{F}_2=\{P: \,\, prime \,\, ideal \,,\ such \,\, that\,\, S [[x;\beta,u]]x\nsubseteq P\}$.\end{center}  Note that $\displaystyle\bigcap_{P\in \mathcal{F}_1}P=Nil_{*}(S)\oplus \displaystyle\sum _{i\geq 1}Sx^i$. Now, for each strongly $\beta$-prime ideal $Q$ of $S$, we have by similar methods of Corollary \ref{primeideals}, that $Q [[x;\beta,u]]$ is prime  and we easily get that  each prime ideal $P$ of $\mathcal{F}_2$ implies that $P\cap S$ is a strongly $\beta$-prime ideal of $S$.  Thus, $\displaystyle\bigcap_{P\in \mathcal{F}_2} P  \supseteq N_{\beta}(S)$. Hence, $Nil_{*}(S [[x;\beta,u]]=(\displaystyle\bigcap_{P\in \mathcal{F}_1} P)\cap (\displaystyle\bigcap_{Q\in \mathcal{F}_2}Q)\supseteq  (Nil_{*}(S) +\displaystyle\sum_{i\geq 1}Sx^i)\cap (N_{\beta}(S)) [[x;\alpha,w]])\supseteq Nil_{*}(S)\cap N_{\beta}(S) \oplus \displaystyle\sum_{i\geq 1}N_{\beta}(S)x^i$.

On the other hand,  since for each prime ideal $L$ of $S$ we have that $L\oplus \displaystyle\sum _{i\geq 1}Sx^i$ is a prime ideal of $S [[x;\beta,u]]$ and  in the same way of Corollary \ref{primeideals}  we have that $N [[x;\beta,u]]$ is prime for  each strongly $\beta$-prime ideal $N$ of $S$, then  $Nil_{*}(S [[x;\beta,u]]) \subseteq (Nil_{*}(S)\cap N_{\beta}(S)) \oplus \displaystyle\sum_{i\geq 1}N_{\beta}(S)x^i$.

So, $Nil_{*}(S [[x;\beta,u]])=Nil_{*}(S)\cap N_{\beta}(S) \oplus \displaystyle\sum_{i\geq 1}N_{\beta}(S)x^i$.

\end{dem}

\begin{prop}\label{primeradical4} Let $\alpha$ be an unital  twisted partial action with enveloping action $(T,\beta,u)$. Then the prime radical  $Nil_{*}(R [[x;\alpha,w]])$ of $R [[x;\alpha,w]]$ is $Nil_{*}(R [[x;\alpha,w]])=(N_{\alpha}(R)\cap Nil_{*}(R))\oplus \displaystyle\sum_{i\geq 1}(N_{\alpha}(R)\cap D_i)x^i$, where $N_{\alpha}(R)$ is the intersection of all strongly $\alpha$-prime ideals of $R$.\end{prop}

\begin{dem}  Using the same methods of  (\cite{wag e fer}, Lemma 2.9),  we have,  for each strongly $\beta$-prime ideal $Q$ of $T$,  that $Q\cap R$ is a strongly $\alpha$-prime ideal of $R$ and since $R$ is an ideal of $T$ we have that $Nil_{*}(T)\cap R=Nil_{*}(R)$. Thus, by   Corollary \ref{primeradical} we easily get that  $Nil_{*}(R [[x;\alpha,w]])=Nil_{*}(T [[x;\beta,u]])\cap R [[x;\alpha,w]]= (Nil_{*}(R)\cap N_{\alpha}(R))\oplus \displaystyle\sum_{i\geq 1}(Nil_{\alpha}(R)\cap D_i)x^i$\end{dem}

\begin{obs} In the last result,  we use the fact that the twisted partial action has an enveloping action, but we do not know  if the Proposition \ref{primeradical4} is true for unital twisted partial actions of $\mathbb{Z}$ without enveloping action. To solve this problem we need to know if the following result is true:

\begin{center} Let $P$ be a prime ideal of $R [[x;\alpha,w]]$ such that $1_ix^i\notin P$ for some $i\geq 1$. Then $P\cap R$ is strongly $\alpha$-prime.\end {center}

\end{obs}

As it happened in the   (\cite{CFMH}, Example 2.6) we obtain by a similar example that  the twisted partial skew power series over semiprime Goldie rings are  not necessary semiprime, but,  if we input the condition  of ``finite type" we get the following.


\begin{teo} Let $\alpha$ be an unital  twisted partial action. If  $R$ is semiprime Goldie and  $\alpha$ is a twisted partial action of finite type,  then $R [[x;\alpha,w]]$ is semiprime Goldie.
\end{teo}

\begin{dem}   Since $\alpha$ is of finite type and  $rank(R)$ is finite, then by Theorem  \eqref{envolvente} $\alpha$ has enveloping action  $\beta=\left(B, \{\beta _g\}_{g\in G}, \{u_{(g,h)}\}_{(g,h)\in {G\times G}}\right)$. In this case, by (\cite{lmsw}, Corollary 4.18),  $T$ is semiprime Goldie and we claim that   $T [[x;\beta,u]]$ is semiprime. In fact, suppose that $Nil_{*}(T [[x;\beta,u]])$  is not zero. Then, by (\cite{lam1}, Lemma 10.10.29), $Nil_{*}(T [[x;\beta,u]])$ contains a nonzero nilpotent ideal $L$, since by Theorem 3.4 we have that $T[[x;\beta,u]]$ is Goldie.  By the fact that $T$ is semiprime we have that $Nil_{*}(T [[x;\beta,u]])=\displaystyle\sum_{i\geq 1}N_{\beta}(T)x^i$. Now, consider $ H=\{0\neq a\in N_{\beta}(T):\exists 0\neq f\in  L\,\, such\,\, that  \,\, f=ax^{j}+...\in L\}\cup {0}$.  It is not difficult to see that $H$ is a nonzero ideal of $T$ with  $ \beta_i(H)\subseteq H$, for all $i\in\mathbb{Z}$. Since $L$ is nilpotent, we obtain that  $H$ is nilpotent and consequently  $H=0$,  because of  $T$ is semiprime, which is a contradiction. So,  $Nil_{*}(T [[x;\beta,u]])=0$.

By Corollary \ref{primeradical} we have that  $$Nil_{*}(T [[x;\beta,u]])\cap R [[x;\alpha,w]]=Nil_{*}(R [[x;\alpha,w]])$$ which implies that    $Nil_{*}(R [[x;\alpha,w]])=0$. Therefore, $R [[x;\alpha,w]]$ is semiprime  Goldie.
\end{dem}

\end{document}